\documentclass[a4paper,11pt]{article}
\usepackage{latexsym}
\usepackage{amsmath}
\usepackage{amssymb}
\usepackage{enumerate}
\usepackage{theorem}
\usepackage{array}
\newtheorem{theorem}{Theorem}[section]
\newtheorem{lemma}[theorem]{Lemma}
\newtheorem{corollary}[theorem]{Corollary}
\newtheorem{proposition}[theorem]{Proposition}
\newtheorem{atheorem}{Theorem A.\!\!}
\newtheorem{alemma}[atheorem]{Lemma A.\!\!}
\theorembodyfont{\rmfamily}
\newtheorem{definition}[theorem]{Definition}
\newtheorem{notation}[theorem]{Notation}
\newtheorem{setup}[theorem]{Setup}

\newtheorem{aremark}[atheorem]{Remark A.\!\!}
\newtheorem{aexample}[atheorem]{Example A.\!\!}
\newtheorem{notationinintro}{Notation}

\newcommand{\proof}{\noindent \mbox{\em Proof.\hspace*{2mm}}}
\newcommand{\qed}{\hfill \mbox{$  \Box $}}
\newcommand{\ssgyokan}{\vskip 4pt}

\title{The characterization of Hermitian surfaces by the number of points}
\author{
Masaaki Homma
\thanks{Partially supported by Grant-in-Aid
for Scientific Research (24540056), JSPS.}
\\
 Department of Mathematics and Physics\\
Kanagawa University\\
Hiratsuka 259-1293, Japan\\
homma@kanagawa-u.ac.jp
\and
Seon Jeong Kim
\thanks{Partially supported by Basic Science Research Program through the National Research Foundation of Korea(NRF) 
funded by the Ministry of Education, Science and Technology (2012R1A1A2042228).
}\\
 Department of Mathematics and RINS\\
Gyeongsang National University\\
Jinju 660-701, Korea \\
skim@gnu.kr
}
\date{}

\begin{document}
\maketitle
\begin{abstract}
The nonsingular Hermitian surface of degree $\sqrt{q} +1$
is characterized by its number of $\Bbb{F}_q$-points among
the surfaces over $\Bbb{F}_q$ of degree $\sqrt{q} +1$
in the projective 3-space
without $\Bbb{F}_q$-plane components.
\end{abstract}

\section{Introduction}
Hermitian varieties are known as ones having particular properties
over finite fields.
Throughout this paper, except Section~3, $q$ is an even power of
a prime number $p$.
A Hermitian variety over ${\Bbb F}_q$ is a hypersurface in ${\Bbb P^n}$
defined by
\begin{equation}\label{def_hermitian}
(X_0^{\sqrt{q}}, \dots , X_n^{\sqrt{q}}) A \, {}^t\!(X_0, \dots , X_n) =0,
\end{equation}
where $A$ is a square matrix of order $n+1$ whose entries are in ${\Bbb F}_q$
with the property ${}^t\!A=A^{(\sqrt{q})}$;
here $A^{(\sqrt{q})}$ means taking entry-wise
the $\sqrt{q}$-th power, and ${}^t\!A$ is the transposed matrix of $A$.
We refer this kind of polynomial as a Hermitian polynomial over ${\Bbb F}_q$.
For a homogeneous polynomial $F$ of degree $\sqrt{q}+1$
whose coefficients are in ${\Bbb F}_q$, the hypersurface given by $F=0$ is
Hermitian if and only if there is an element $\rho \in {\Bbb F}_q^{\ast}$
such that $\rho F$ is a Hermitian polynomial.
The family of Hermitian polynomials over ${\Bbb F}_q$
forms an ${\Bbb F}_{\! \sqrt{q}}$ vector space.
It is obvious that the Hermitian polynomial (\ref{def_hermitian})
defines a nonsingular Hermitian variety if and only if $\det A \neq 0$.
By the standard argument,
a nonsingular Hermitian variety is projectively equivalent to
the variety
\begin{equation}\label{standard_form}
H_{n-1} : X_0^{\sqrt{q}+1} + X_1^{\sqrt{q}+1} 
+ \dots  + X_n^{\sqrt{q}+1} =0
\end{equation}
over ${\Bbb F}_q$.

The number $N_q(H_{n-1})$ of ${\Bbb F}_q$-points of $H_{n-1}$ is
\begin{equation}\label{number_points_hermitian}
N_q(H_{n-1}) = (\sqrt{q}^{n+1} + (-1)^{n})(\sqrt{q}^{n} + (-1)^{n-1})/(q-1),
\end{equation}
which is due to \cite{bos-cha1966}.

This number is remarkable in the following sense.
The Weil conjecture \cite{wei1949} established by Deligne \cite{del1973}
implies that the number $N$ of ${\Bbb F}_q$-points of
a nonsingular hypersurface of degree $d$ in ${\Bbb P}^n$
defined over ${\Bbb F}_q$ is bounded by
\begin{equation}\label{weil_deligne}
N \leq \frac{q^n -1}{q-1} + 
 \frac{d-1}{d} \left( (d-1)^n - (-1)^n \right) \sqrt{q}^{n-1}.
\end{equation}
Furthermore, if equality holds in (\ref{weil_deligne})
for a certain nonsingular hypersurface of degree $d$,
then $d \leq \sqrt{q}+1$.

As for this additional claim, see \cite[Corollary 2.2]{sor1994}
or \cite[Corollary 4.3]{hom-kim2013}.
The number $N_q(H_{n-1})$
achieves the equality in (\ref{weil_deligne})
for $d=\sqrt{q} +1$.

There is a characterization of the Hermitian curve
among the family of plane curves over ${\Bbb F}_q$ of degree $\sqrt{q} +1$
without ${\Bbb F}_q$-line components \cite{hir-sto-tha-vol}.
Note that $N_q(H_{1}) = \sqrt{q}^3 +1$.

\begin{theorem}[Hirschfeld-Storme-Thas-Voloch]\label{hermitiancurve}
Suppose that $q \neq 4$.
Let $C$ be a plane curve over ${\Bbb F}_q$ of degree $\sqrt{q} +1$
without ${\Bbb F}_q$-line components.
If $C$ has $\sqrt{q}^3 +1$ points over ${\Bbb F}_q$,
then $C$ is a Hermitian curve.
\end{theorem}

The purpose of this paper is to show a similar fact
for surfaces in ${\Bbb P}^3$.

\begin{theorem}\label{maintheorem}
Let $S$ be a surface in ${\Bbb P}^3$
defined over ${\Bbb F}_q$ without ${\Bbb F}_q$-plane components.
If the degree of $S$ is $\sqrt{q}+1$ and
$N_q(S) = (\sqrt{q}^3+1)(q+1)$,
then $S$ is a nonsingular Hermitian surface over ${\Bbb F}_q$.
\end{theorem}

After Hirschfeld-Storme-Thas-Voloch's characterization of the Hermitian curve
established, R\"{u}ck-Stichtenoth gave another characterization 
of the Hermitian curve among
the family of nonsingular curves defined over ${\Bbb F}_q$ of genus
$\frac{\sqrt{q}(\sqrt{q}-1)}{2}$ \cite{ruc-sti1994}.
A connection of those two characterization will be mentioned in Appendix.

\begin{notationinintro}
\begin{itemize}
\item When $X$ is an algebraic set in ${\Bbb P}^n$
defined by equations over ${\Bbb F}_q$,
the set of ${\Bbb F}_q$-points in $X$ is denoted by $X({\Bbb F}_q)$,
and the cardinality of $X({\Bbb F}_q)$ by $N_q(X)$.
\item When coordinates $X_0, \dots , X_n$ of ${\Bbb P}^n$ are given,
for a homogeneous polynomial $h$ in those variables,
$\{ h=0 \}$ means the hypersurface defined by $h=0$.
\item Let $D$ and $E$ are curves in ${\Bbb P}^2$ without common components.
For a point $Q \in D\cap E$,
$i(D.E;Q)$ denotes the local intersection multiplicity of $D$ and $E$ at $Q$,
and $(D.E) = \sum_{Q \in D \cap E} i(D.E;Q)$.
\item For a plane curve $C \subset {\Bbb P}^2$
and a nonsingular point $P \in C$,
$T_P(C)$ denotes the embedded tangent line to $C$ at $P$.
\item The multiplicative set ${\Bbb F}_q \setminus \{0\}$ is denoted by
${\Bbb F}_q^{\ast}$.
\item For a finite set $Y$, ${}^{\#}Y$ denotes the number of elements of $Y$.
\end{itemize}
\end{notationinintro}
\section{Plane curves}
For a plane curve $C$ over ${\Bbb F}_q$,
we proved a simple bound for $N_q(C)$
in a series of papers \cite{hom-kim2009, hom-kim2010a, hom-kim2010b},
which had been originally conjectured by Sziklai~\cite{szi2008}.
\begin{theorem}[Sziklai bound]\label{sziklai_bound}
Let $d$ be an integer with $2 \leq d \leq q+2$,
and $C$ a curve of degree $d$ in ${\Bbb P}^2$
defined over ${\Bbb F}_q$ without ${\Bbb F}_q$-linear components.
Then
\begin{equation}\label{ineq_sziklai_bound}
N_q(C) \leq (d-1)q+1
\end{equation}
except for curves over ${\Bbb F}_4$
which are projectively equivalent to the curve defined by
\[
K: \ 
(X_0 + X_1 + X_2)^4 + (X_0X_1 + X_1X_2 +X_2X_0)^2 +X_0X_1X_2(X_0 + X_1 + X_2)
=0.
\]

For the exceptional curve $K$ above, $N_4(K)=14$.
\end{theorem}

The number $N_q(H_1)$ attains the equality of the Sziklai bound too.

This bound implies a sufficient condition on a plane curve over ${\Bbb F}_q$
to be absolutely irreducible.

\begin{corollary}\label{absolutelyirreducible}
Let $C$ be a curve over ${\Bbb F}_q$ of degree $d$ in ${\Bbb P}^2$
without ${\Bbb F}_q$-linear components.
If $N_q(C) \geq (d-2)q +3$, then $C$ is absolutely irreducible.
\end{corollary}
\proof
First note that the range of $d$ is $2 \leq d \leq q+2$ by assumptions.

(Step~1) Suppose $C$ is decomposed into two curves $C_1$ and $C_2$ over
${\Bbb F}_q$.
For each $i =1, 2$, let $d_i$ be the degree of $C_i$.
If either ``$q\neq 4$" or 
``$q=4$ and none of the $C_i$'s is projectively equivalent to $K$,"
then
\[
N_q(C) \leq N_q(C_1)+N_q(C_2)
  \leq ((d_1-1)q + 1) +((d_2-1)q + 1)
\leq (d-2)q + 2.
\]

We consider the exceptional case, namely
``$q=4$ and one of the $C_i$'s is projectively equivalent to $K$
oner ${\Bbb F}_4$."
Applying an ${\Bbb F}_4$-linear transformation if necessary,
we may assume that $C_1=K$.
Then $d=6$ and $d_2=2$, because $2 \leq d \leq 6$ and
$C$ has no ${\Bbb F}_4$-linear components.
Since $C_2$ is ${\Bbb F}_4$-irreducible of degree $2$,
$N_4(C_2) =5$ or $1$.
When $N_4(C_2) =1$, $N_4(C) \leq 14 +1 < (6-2)4 + 3.$
If $C_2$ is an absolutely irreducible conic,
then ${}^{\#}(C_2 \cap \mathbb{P}^2(\mathbb{F}_2)) \leq 4$,
because the maximum length of $\mathbb{F}_2$-arc is $4$.
Since $K(\mathbb{F}_4) =
\mathbb{P}^2(\mathbb{F}_4) \setminus \mathbb{P}^2(\mathbb{F}_2)$
(see \cite[\S 3]{hom-kim2009}),
${}^{\#}(K(\mathbb{F}_4) \cap C_2(\mathbb{F}_4))\geq 1$.
Hence
$N_q(C) \leq 14 + 5 -1 = (6-2)4+2.$

Therefore $C$ is irreducible over ${\Bbb F}_q$.

(Step~2)
Suppose that the ${\Bbb F}_q$-irreducible curve $C$
is not absolutely irreducible.
Let $D$ be an irreducible component of $C$.
Then $D$ is defined over an extension $\mathbb{F}_{q^t}$ of $\mathbb{F}_q$
with $t \geq 2$, and
$C= D \cup D^{(q)} \cup \dots \cup D^{(q^{t-1})}$,
where $D^{(q)},  \dots , D^{(q^{t-1})}$
are conjugates of $D$ over $\mathbb{F}_q$.
Hence $C({\Bbb F}_q) \subset D \cap D^{(q)}$ and
$\deg D = \deg D^{(q)} \leq \frac{d}{2}.$
Hence $N_q(C) \leq \left( \frac{d}{2}\right)^2$.
Since $2 \leq d \leq q+2$,
we have $\left( \frac{d}{2}\right)^2< (d-2)q + 3$.
Indeed, consider the quadratic
$P(d) = \left( \frac{d}{2}\right)^2 - (d-2)q - 3.$
Then $P(2)<0$ and $P(q+2)<0$, which implies
$P(d)<0$ for any $d$ with $2 \leq d \leq q+2$.

Therefore the curve $C$ is absolutely irreducible too.
\qed

\ssgyokan

An important step of the proof of Theorem~\ref{maintheorem} is
to study plane sections of $S$
and to apply the Hirschfeld-Storme-Thas-Voloch theorem to some of them,
but the exception in the assertion of this theorem harms this step for
$q=4$.
So we take a detour in this case.

\begin{lemma}\label{detour}
Let $C$ be a plane curve of degree $3$ over $\mathbb{F}_4$
without $\mathbb{F}_4$-linear components.
Suppose $N_4(C)=9$.
Then $C$ is absolutely irreducible
and all $\mathbb{F}_4$-points are nonsingular.
In addition, suppose that all $\mathbb{F}_4$-points are flexes,
i.e., $i(T_p(C).C;P) =3$ for any $P \in C(\mathbb{F}_4)$.
Then $C$ is a Hermitian curve.
\end{lemma}
\proof
From (\ref{absolutelyirreducible}),
$C$ is absolutely irreducible.
Let $P \in C(\mathbb{F}_4)$.
Assume that $P$ is a singular point.
Then for any $\mathbb{F}_4$-line passing through $P$,
${}^{\#}(l(\mathbb{F}_4)\cap C \setminus \{P\}) \leq 1$
because $(l.C)=3$ and $i(l.C; P) \geq 2$.
Hence $N_4(C) \leq 6$ because
$C(\mathbb{F}_4) = \cup_{l \in \Check{P}} (l(\mathbb{F}_4)\cap C)$,
which contradicts to the assumption $N_4(C)=9$.
Hence all $\mathbb{F}_4$-points are nonsingular.

Next we show the additional statement.
Since $\deg C=3$ and $i(T_p(C).C;P) =3$ for any $P \in C(\mathbb{F}_4)$,
a line joining two $\mathbb{F}_4$-points of $C$ meets with $C$
at the third $\mathbb{F}_4$-point.
Hence the $9$ $\mathbb{F}_4$-points $C(\mathbb{F}_4)$
can be divided three triples
\[
\{P_{01}, P_{02}, P_{03}\} \ \{P_{11}, P_{12}, P_{13}\} \ 
\{P_{21}, P_{22}, P_{23}\}
\]
such that three points of each triple are collinear.
Let $l_i$ be the line on which $P_{i1}, P_{i2}$ and  $P_{i3}$ lie, and
$l_i(\mathbb{F}_4) = \{P_{i1}, P_{i2}, P_{i3}, Q_i, R_i \}$.
Assume that $l_0, l_1$ and $l_2$ are concurrent.
We may assume that $Q_0=Q_1=Q_2$.
Hence $T_{P_{0 \mu}} \cap l_1$ must be $R_1$ for any $\mu =1,2,3$, and also
$T_{P_{2 \nu}} \cap l_1 = R_1$ for any $\nu =1,2,3$
because $P_{i \alpha}$ $(i, \alpha =1,2,3)$ are flexes.
But, since the number of $\mathbb{F}_4$-lines passing through $R_1$ is $5$
and $T_{P_{i \mu}} \neq T_{P_{j \nu}} $ for any pair $(P_{i \mu}, P_{j \nu})$,
it is impossible.

Therefore
$\{l_0, l_1, l_2\}$ forms a frame of $\mathbb{P}^2$,
that is, we may assume that
$l_i =\{ X_i=0 \}$, where $X_0, X_1, X_2$ are coordinates of $\mathbb{P}^2$.
Let $C=\{ F(X_0, X_1, X_2)=0\}$.
Since
\begin{gather}
 F(X_0, X_1, 0) = c_2(X_0^3 + X_1^3) \ \ 
            (c_2 \in \mathbb{F}_4^{\ast}) \notag \\
 F(X_0, 0, X_2) = c_1(X_0^3 + X_2^3) \ \ 
            (c_1 \in \mathbb{F}_4^{\ast}) \notag \\
 F(0, X_1, X_2) = c_0(X_1^3 + X_2^3) \ \ 
            (c_0 \in \mathbb{F}_4^{\ast}), \notag 
\end{gather}
we have
 $F(X_0, X_1, X_2) = c(X_0^3 + X_1^3+X_2^3) + \gamma X_0 X_1 X_2$
with $c \in \mathbb{F}_4^{\ast}$ and $\gamma \in \mathbb{F}_4$.
If $\gamma\neq 0$, then $(1,1,\frac{c}{\gamma})$
is also an $\mathbb{F}_4$-point of $C$,
which is not any of the $P_{i, \alpha}$'s.
This is a contradiction.
Hence $\gamma =0$, that is, $C$ is Hermitian.
\qed

We close this section with a remark on plane curves
that are used in the proof of the main theorem.

\begin{lemma}[Segre]\label{segre_bound}
Let $d$ be an integer with $1 \leq d \leq q+1$,
and $C$ be a curve of degree $d$ in ${\Bbb P}^2$
defined over ${\Bbb F}_q$,
which may have ${\Bbb F}_q$-linear components.
Then $N_q(C) \leq dq+1$,
and equality holds if and only if $C$ is a pencil of
$d$ ${\Bbb F}_q$-lines.
\end{lemma}
\proof
See Segre \cite[II, \S 6 Observation IV ]{seg1959}
or Homma-Kim \cite[Remark 1.2]{hom-kim2009}.
\qed

\section{An elementary bound}
In this section, $q$ is simply a power of $p$, that is,
it need not be an even power of $p$.

In \cite{hom-kim2013},
we established an upper bound for $N_q(X)$ of a hypersurface $X$
over ${\Bbb F}_q$ without ${\Bbb F}_q$-linear components,
particularly a surface in ${\Bbb P}^3$ without
${\Bbb F}_q$-plane components.

\begin{theorem}
Let $S$ be a surface of degree $d$ in ${\Bbb P}^3$
defined over ${\Bbb F}_q$
without ${\Bbb F}_q$-plane components.
Then
\begin{equation}\label{elementary_bound}
N_q(S) \leq (d-1)q^2 + dq +1.
\end{equation}
\end{theorem}

We refer the bound (\ref{elementary_bound})
as the elementary bound.
When $d=\sqrt{q}+1$,
the Weil-Deligne bound (\ref{weil_deligne}) for $n=3$ agrees with
the elementary bound.
In this section,
we investigate the geometry of a surface in ${\Bbb P}^3$
whose number of ${\Bbb F}_q$-points achieves the bound
(\ref{elementary_bound}).
There are at least two examples other than the Hermitian surface
each of which attains the equality in (\ref{elementary_bound})
\footnote{
After wrote up the first draft of this paper,
we proved that these two examples together with
the Hermitian surface are only examples having this property
\cite{hom-kim2015a},
in which we referred this paper.
Recently, Tironi \cite{tir2014preprint} has settled
the classification of hypersurfaces in $\mathbb{P}^n$
that achieve the elementary bound stated in \cite{hom-kim2013}.
}.

From now on, we keep the following setup
until the end of this section.

\begin{setup}
Let $d$ be an integer with $2\leq d \leq q+1$.
Let $S$ be a surface of degree $d$ in ${\Bbb P}^3$
defined over ${\Bbb F}_q$
without ${\Bbb F}_q$-plane components.
Furthermore we suppose that $N_q(S)$
achieves the equality in (\ref{elementary_bound}).
\end{setup}
\begin{lemma}\label{contains_line}
The surface $S$ contains an ${\Bbb F}_q$-line.
\end{lemma}
\proof
Suppose $S$ does not contain any ${\Bbb F}_q$-lines.
Let $H$ be any ${\Bbb F}_q$-plane in ${\Bbb P}^3$.
Then $S\cap H$ is a plane curve of degree $d$ over ${\Bbb F}_q$
in $H = {\Bbb P}^2$, and has no ${\Bbb F}_q$-line as a component.
Hence
\begin{equation}
 N_q(S\cap H) \leq
  \left\{
      \begin{array}{cl}
       (d-1)q +1 & \mbox{\rm if $(d,q) \neq (4,4)$}\\
       14        & \mbox{\rm if $(d,q) = (4,4)$}
      \end{array}
  \right.
\end{equation}
by Lemma~\ref{sziklai_bound}.
In a term of \cite{hom2012},
defining the $s$-degree $\delta$ of $S({\Bbb F}_q)$ by
\[
\delta = \max \{ N_q(S\cap H) \mid \mbox{\rm $H$ is an ${\Bbb F}_q$-plane} \}, 
\]
we have
\[
N_q(S) \leq (\delta -1)q + 1 +
\left\lfloor
\frac{\delta -1}{q +1}
     \right\rfloor
\]
by \cite[Proposition 2.2]{hom2012}.
Hence if $(d,q) \neq (4,4)$, then
\begin{eqnarray*}
 N_q(S) & \leq & (d-1)q^2 +1 + \left\lfloor
\frac{(d -1)q}{q +1}
     \right\rfloor\\
     &=& (d-1)q^2 +1 + \left\lfloor
\frac{(d -2)(q+1)+q+2-d}{q +1}
     \right\rfloor\\
     &=& (d-1)(q^2+1) \, {\rm ;}
\end{eqnarray*}
and if $(d,q) = (4,4)$, then $N_q(S) \leq 55$.
In either case, $N_q(S)$ can't be
$(d-1)q^2 + dq +1.$
\qed

\begin{definition}
Let $l_1, \dots , l_d$ be ${\Bbb F}_q$-lines in ${\Bbb P}^3$
with $d \geq 2$.
The union of those $d$-lines $Z = \cup_{i=1}^d l_i$
is called a planar ${\Bbb F}_q$-pencil of degree $d$
if those $d$-lines lie on a plane simultaneously
and $l_1 \cap \dots \cap l_d \neq \emptyset$.
The ${\Bbb F}_q$-point $\{v_Z\} = l_1 \cap \dots \cap l_d$
is called the vertex of $Z$.
\end{definition}
\begin{notation}
For an ${\Bbb F}_q$-line $l$ of ${\Bbb P}^3$,
the set of ${\Bbb F}_q$-planes containing the line $l$
is denoted by $\check{l}({\Bbb F}_q)$.
\end{notation}
\begin{lemma}\label{planar_pencil}
Let $l$ be an ${\Bbb F}_q$-line on the surface $S$.
\begin{enumerate}[{\rm (i)}]
\item If an ${\Bbb F}_q$-plane $H$ contains $l$,
then $S\cap H$ is a planar ${\Bbb F}_q$-pencil of degree $d$.
\item The map
$\check{l}({\Bbb F}_q) \ni H \mapsto v_{S\cap H} \in l({\Bbb F}_q)$
is bijective.
\end{enumerate}
\end{lemma}
\proof
(i) Since $S$ has no ${\Bbb F}_q$-plane components,
$S\cap H$ is a plane curve of degree $d$ in $H$,
and $N_q(S \cap H) \leq dq+1$
by Lemma~\ref{segre_bound}.
Counting the cardinality of $S({\Bbb F}_q)$ by the decomposition
\[
S({\Bbb F}_q) = \cup_{H \in \check{l}({\Bbb F}_q)}
   \left(
     (S\cap H)({\Bbb F}_q) \setminus l({\Bbb F}_q)
   \right)
   \bigcup l({\Bbb F}_q),
\]
we have
\begin{eqnarray*}
(d-1)q^2 + dq +1=N_q(S) &\leq& (q+1)(dq+1-(q+1))+(q+1)\\
                        &=& (d-1)q^2 + dq +1.
\end{eqnarray*}
Hence $N_q(S\cap H) =dq +1$.
So $S \cap H$ is a planar ${\Bbb F}_q$-pencil of degree $d$
by the latter part of Lemma~\ref{segre_bound}.

(ii) Let
$\check{l}({\Bbb F}_q)
=\{ H_1, H_2, \dots , H_{q+1} \}.$
Since $S\cap H_i$ is a planar  ${\Bbb F}_q$-pencil and
has the line $l$ as a component,
we may set notations as
\[
S \cap H_i = l \cup l_{i,1} \cup \dots \cup l_{i, d-1}
\]
and
\[ v_i = v_{S \cap H_i} .\]
It is obvious that $v_i \in l({\Bbb F}_q).$
Since $ \check{l}({\Bbb F}_q) $ and $l({\Bbb F}_q)$
have the same cardinality,
it is enough to show that the map
$\check{l}({\Bbb F}_q) \ni H_i \mapsto v_i \in l({\Bbb F}_q)$
is surjective.

Contrary, suppose this map is not surjective.
Pick a point $Q \in l({\Bbb F}_q) \setminus \{ v_1, \dots , v_{q+1} \},$
and choose an ${\Bbb F}_q$-plane $K$ such that $K \ni Q$
and $ K \not\supset l$.
Then we face two consequences:
\begin{enumerate}
\item[($\alpha$)] $S \cap K$ does not contain any ${\Bbb F}_q$-line;
\item[($\beta$)] $N_q(S \cap K) \geq (d-1)q +d$,
\end{enumerate}
as we verify them below.
Since $S \cap K$ is a plane curve of degree $d$ in $K$,
those two conditions are incompatible each other
by Lemma~\ref{sziklai_bound}.

{\em The verification of} ($\alpha$).
If $S \cap K$ contains an ${\Bbb F}_q$-line,
then it must be a planar ${\Bbb F}_q$-pencil by (i).
Since $Q$ is a point of $S$, there is an ${\Bbb F}_q$-line
$m$ passing through $Q$ among $d$ lines of $S \cap K$.
Hence $l$ and $m$ spans an ${\Bbb F}_q$-plane
which is one of the $H_i$'s, say $H_i$.
Then $v_i \in l \cap m =\{ Q \}$,
which contradicts to the choice of $Q$.

{\em The verification of} ($\beta$).
Let $Q_{i,j}$ be the intersection point of $l_{i,j}$ with $K$.
Then the plane containing $l$ and $Q_{i,j}$ is $H_i$
and the line containing $v_i$ and $Q_{i,j}$ is $l_{i,j}$.
Hence $Q_{i,j} =Q_{i',j'}$
implies $(i,j) = (i',j')$.
Since
\[
(S\cap K)({\Bbb F}_q) \supset 
  \{Q_{i,j} \mid 1 \leq i \leq q+1, \, 1 \leq j \leq d-1 \} \cup \{Q\},
\]
we have $N_q(S\cap K) \geq (d-1)(q+1) +1.$
\qed

\begin{lemma}
Let $H$ be an ${\Bbb F}_q$-plane
such that $S \cap H$ is a planar ${\Bbb F}_q$-pencil of degree $d$.
If an ${\Bbb F}_q$-line $l \subset S$ passes through $v_{S\cap H}$,
then $l$ is a component of $S \cap H$.
\end{lemma}
\proof
Since $d \geq 2$,
there are two distinct components $l_1$ and $l_2$ of $S \cap H$.
Suppose that $l$ is not contained in $H$.
Then  $l$ and $l_i$ span an ${\Bbb F}_q$-plane, say $H_i$.
Hence $S\cap H_i$ is also a planar ${\Bbb F}_q$-pencil of degree $d$
by Lemma~\ref{planar_pencil}~(i).
By the construction of $H_1$ and $H_2$,
\[
v_{S \cap H_1} = l \cap l_1 = v_{S \cap H} = l \cap l_2 = v_{S \cap H_2},
\]
which contradicts (ii) of Lemma~\ref{planar_pencil}.
Hence $l$ is contained in $H$, and hence it is a component of $S\cap H$.
\qed
\begin{corollary}\label{bijective_pencil_point}
For any ${\Bbb F}_q$-point $P$ of $S$,
there is a unique ${\Bbb F}_q$-plane $H$
such that $S\cap H$ is a planar ${\Bbb F}_q$-pencil
of degree $d$ with $v_{S\cap H} =P$.
\end{corollary}
\proof
From Lemma~\ref{contains_line} and Lemma~\ref{planar_pencil} (i),
there is a line $l$ on $S$ containing $P$.
Indeed the existence of an ${\Bbb F}_q$-line $l'$ on $S$
is guaranteed by Lemma~\ref{contains_line};
if the point $P$ in question lies on $l'$, there is nothing to do;
when $P \not\in l'$,
consider the ${\Bbb F}_q$-plane $H'$ containing $l'$ and $P$,
then $S \cap H'$ is a planar ${\Bbb F}_q$-pencil
by Lemma~\ref{planar_pencil} (i), which contains $P$;
hence there is a line $l$ on $S$ containing $P$.
Hence, by (ii) of Lemma~\ref{planar_pencil}
there is a desired ${\Bbb F}_q$-plane,
and such a plane is uniquely determined by $P$.
\qed

\section{A characterization of Hermitian surfaces}
In this section,
$q$ is assumed to be an even power of a prime number $p$ again,
and $S$ is a surface in ${\Bbb P}^3$ of degree $\sqrt{q}+1$
defined over ${\Bbb F}_q$ without ${\Bbb F}_q$-plane
components such that $N_q(S) = (\sqrt{q}^{3}+1)(q+1)$.
Note that this number $N_q(S)$
achieves both the Weil-Deligne bound (\ref{weil_deligne})
for $n=3$ and the elementary bound (\ref{elementary_bound}).

\begin{proposition}\label{sqrt_q_plus_one_case}
With the above situation,
let $H$ be an ${\Bbb F}_q$-plane of ${\Bbb P}^3$.
Then either
\begin{enumerate}[{\rm (1)}]
\item $S \cap H$ is a planar ${\Bbb F}_q$-pencil of degree $\sqrt{q}+1$, or
\item $S \cap H$ is a nonsingular Hermitian curve
of degree $\sqrt{q}+1$ over ${\Bbb F}_q$.
\end{enumerate}
Furthermore, let $\nu_1$ denote the number of
${\Bbb F}_q$-planes having the property {\rm (1)} above,
and $\nu_2$ the property {\rm (2)}.
Then
\begin{equation*}
\left\{
\begin{array}{ccc}
 \nu_1 &=& N_q(S) = (\sqrt{q}^3 +1)(q+1) \\
 \nu_2 &=& \sqrt{q}^3 (q+1)(\sqrt{q} -1).
\end{array}
\right.
\end{equation*}
\end{proposition}
\proof
Let $\check{\Bbb P}^3({\Bbb F}_q)$ denote
the set of ${\Bbb F}_q$-planes of ${\Bbb P}^3$,
and $\mathcal{P}_1$ the set of ${\Bbb F}_q$-planes
having the property (1).
The map $\mathcal{P}_1 \ni H \mapsto v_{S \cap H} \in S({\Bbb F}_q)$
is bijective by Corollary~\ref{bijective_pencil_point}.
Hence $\nu_1 = N_q(S)$.
Let $\mathcal{P}_2= \check{\Bbb P}^3({\Bbb F}_q) \setminus \mathcal{P}_1 $.
Note that if $H \in \mathcal{P}_2$, then
$S \cap H$ has no ${\Bbb F}_q$-liner components by
Lemma~\ref{planar_pencil}~(i).
Consider the correspondence
\[
\mathcal{A} = \{(P, H) \mid P \in H\}
\subset S({\Bbb F}_q) \times \check{\Bbb P}^3({\Bbb F}_q)
\]
with two projections
$\pi_1: \mathcal{A} \to  S({\Bbb F}_q)$ and
$\pi_2: \mathcal{A} \to  \check{\Bbb P}^3({\Bbb F}_q)$.
Counting the cardinality of $\mathcal{A}$ by $\pi_1$, we have
\[
{}^{\#}\mathcal{A} = N_q(S)\cdot(q^2+q+1).
\]
If $H \in \mathcal{P}_1$,
then 
$
{}^{\#}\pi_2^{-1}(H) = N_q(S \cap H) = \sqrt{q}^3+q +1,
$
and if $H \in  \mathcal{P}_2$, then
$
{}^{\#}\pi_2^{-1}(H) = N_q(S \cap H) \leq \sqrt{q}^3 +1
$
by (\ref{sziklai_bound}).
Hence, counting the cardinality of $\mathcal{A}$ by $\pi_2$, we have
\begin{eqnarray*}
{}^{\#}\mathcal{A} &\leq & N_q(S)\cdot (\sqrt{q}^3+q +1)
    + \left(
     \frac{q^4 -1}{q-1} -N_q(S)
    \right) (\sqrt{q}^3 +1)\\
    &=& N_q(S)\cdot(q^2+q+1).
\end{eqnarray*}
Therefore, if $H \in \mathcal{P}_2$,
then $N_q(S \cap H) = \sqrt{q}^3 +1$.
Furthermore, when $H \in \mathcal{P}_2$,
since $S \cap H$ has no ${\Bbb F}_q$-line as a component,
$S\cap H$ is a nonsingular Hermitian curve
by (\ref{hermitiancurve}) if $q \neq 4$. 

We need a little more argument if $q=4$.
We want to apply Lemma~\ref{detour}.
If an $\mathbb{F}_4$-point is not a flex of $S\cap H$,
i.e., $i(T_P(S\cap H). S\cap H; P)=2$,
then $T_P(S\cap H)$ meets with $S\cap H$
at exactly two  $\mathbb{F}_4$-points, because
$\deg S\cap H=3$.
Consider all $\mathbb{F}_4$-planes $\{ H_{\lambda} \}$
in $\mathbb{P}^3$ such that $ H_{\lambda}  \supset T_P(S\cap H_{\lambda}). $
Let $a = {}^{\#}\{ \lambda \mid H_{\lambda} \in \mathcal{P}_1 \}$.
Then ${}^{\#}\{ \lambda \mid H_{\lambda} \in \mathcal{P}_2 \}= 5-a$.
Since $N_4(S \cap H_{\lambda} )$ is either $13$ or $9$,
we have
$ 45 = (13-2)a +(9-2)(5-a) +2$.
Hence there is an $\mathbb{F}_4$-plane $ H_{\lambda} $
such that $S \cap H_{\lambda}$
is a planar $\mathbb{F}_4$-pencil of degree $3$,
and the line $T_{P}(S \cap H)$ lies on $H_{\lambda}$.
But it is impossible, because any $\mathbb{F}_4$-line
on $H_{\lambda}$ meets with the planar pencil at either
$1$, or $3$, or $5$ points.
\qed

\ssgyokan


\noindent 
\mbox{\em Last Step of the Proof of Theorem~\ref{maintheorem}.\hspace*{2mm}}
From Proposition~\ref{sqrt_q_plus_one_case},
there is an ${\Bbb F}_q$-plane $H_{\infty} \subset {\Bbb P}^3$
such that $S \cap H_{\infty}$ is a nonsingular
Hermitian curve.
Choose a system of homogeneous coordinates
$X_0, X_1, X_2, X_3$ of ${\Bbb P}^3$
such that
\begin{enumerate}[(i)]
\item $H_{\infty}$ is given by $X_0=0$; and
\item the plane curve $S \cap H_{\infty}$ in $H_{\infty}={\Bbb P}^2$
is given by
\[
\bar{X}_1 \bar{X}_2^{\sqrt{q}} + \bar{X}_1^{\sqrt{q}} \bar{X}_2
 + \bar{X}_3^{\sqrt{q}+1} =0,
\]
where $\bar{X}_1, \bar{X}_2, \bar{X}_3$ are coordinates
on $H_{\infty}$ induced by $X_1, X_2, X_3$ respectively.
\end{enumerate}

Let $P_1=(0,0,1,0)$ and $P_2=(0,1,0,0)$, both of which are points
on $S \cap H_{\infty}$.
For $\alpha = 1, 2$, the tangent line $L_{\alpha}$
at $P_{\alpha}$ to $S\cap H_{\infty}$ in $H_{\infty}={\Bbb P}^2$
is given by $\bar{X}_{\alpha}=0$.
It is easy to see that $S \cap L_{\alpha}= \{ P_{\alpha} \}$.
Hence for an ${\Bbb F}_q$-plane $H \supset L_{\alpha}$,
\[
{}^{\#} \left(
 S({\Bbb F}_q) \cap (H \setminus L_{\alpha})
\right)
=
\left\{
  \begin{array}{ccl}
  \sqrt{q}^3 + 1 -1 & & \mbox{\rm if $S\cap H$ is Hermitian}\\
  \sqrt{q}^3 +q +1 -1& & \mbox{\rm if $S\cap H$ is an ${\Bbb F}_q$-pencil}.
  \end{array}
\right. 
\]
Counting the number $N_q(S) -1$ by using all ${\Bbb F}_q$-planes
containing $L_{\alpha}$,
we know that there is a unique plane $H_{\alpha, 0} \supset L_{\alpha}$
such that $S \cap H_{\alpha, 0}$ is a planar ${\Bbb F}_q$-pencil
of degree $\sqrt{q}+1$;
and $S \cap H_{\alpha, \lambda }$ is a nonsingular Hermitian curve
for other plane $H_{\alpha, \lambda} \supset  L_{\alpha}$.
By changing coordinates of type
\[
\left\{
  \begin{array}{ccc}
   X_1 & \mapsto & X_1 + a X_0 \\
   X_2 & \mapsto & X_2 + b X_0
  \end{array}
\right.
\]
if necessary,
we may suppose that $H_{\alpha, 0}$
is defined by $X_{\alpha}=0$ for $\alpha =1$ and $2$ respectively.
But the situation on $H_{\infty}$ never change.

Summing up, $S$ is defined by
\[
F(X_0, \dots , X_3)
 = X_0 f(X_0, \dots , X_3) + h(X_1, X_2, X_3) =0
\]
where
\[
h(X_1, X_2, X_3)= X_1 X_2^{\sqrt{q}} + X_1^{\sqrt{q}}X_2 + X_3^{\sqrt{q}+1}
\]
and $\deg f = \sqrt{q}$.

Since $S \cap H_{1, 0} = S \cap \{ X_1 = 0 \}$ is a planar ${\Bbb F}_q$-pencil with
the vertex $P_1= (0,0,1,0)$,
\[
F(X_0, 0, X_2 , X_3) =
    X_0 f(X_0, 0, X_2 , X_3) + X_3^{\sqrt{q}+1}
\]
does not contain $X_2$,
because this polynomial must have the form
$ c \prod_j (X_0+\gamma_j X_3),$
with $c \in {\Bbb F}_q^{\ast}.$
Hence, in $F(X_0, X_1, X_2, X_3)$, any monomial containing $X_2$
also contains $X_1$.
By the same argument on $S \cap H_{2, 0} =S \cap \{ X_2 = 0 \}$,
any monomial containing $X_1$
also contains $X_2$.
Therefore $f(X_0, X_1, X_2, X_3)$
is written as
\[
f(X_0, X_1, X_2, X_3) = g_1(X_0, X_3) + g_2(X_0, \dots , X_3)X_1X_2,
\]
where $\deg g_1 = \sqrt{q}$ and $\deg g_2 = \sqrt{q}-2$.

For the plane $H_{1, \lambda} = \{X_0 = \lambda X_1 \} \ 
(\lambda \in {\Bbb F}_q^{\ast}),$
since $S \cap \{X_0 = \lambda X_1 \}$ is a Hermitian curve,
\begin{eqnarray}
\lefteqn{ \rho F(\lambda X_1, X_1, X_2, X_3)=}\nonumber \\
&& \rho \left(
\lambda X_1 (g_1(\lambda X_1, X_3) + g_2(\lambda X_1, X_1, X_2, X_3)X_1 X_2)
  + h(X_1, X_2, X_3)
\right)
\label{polynomialinproblem}
\end{eqnarray}
is a Hermitian polynomial for some $\rho \in {\Bbb F}_q^{\ast}$.
Since $X_3^{\sqrt{q}+1}$ appears only in $h(X_1, X_2, X_3)$,
the constant $\rho$ must be an element of ${\Bbb F}_{\sqrt{q}}^{\ast}$.
Hence (\ref{polynomialinproblem}) is a Hermitian polynomial
even if $\rho =1$.
Since $h$ itself Hermitian,
\[
\lambda g_1(\lambda X_1, X_3)X_1 + 
      \lambda g_2(\lambda X_1, X_1, X_2, X_3)X_1^2 X_2
\]
must be Hermitian.

Any monomial in $g_2(\lambda X_1, X_1, X_2, X_3)X_1^2 X_2$
never appear in $g_1(\lambda X_1, X_3)X_1$.
So $g_2(\lambda X_1, X_1, X_2, X_3))X_1^2 X_2$
contains the monomials only of types $X_i^{\sqrt{q}+1}$ or $X_i^{\sqrt{q}}X_j$.
Hence
\[
 g_2(\lambda X_1, X_1, X_2, X_3)X_1^2 X_2 = \mu X_1^{\sqrt{q}} X_2
  \ (\mu \in {\Bbb F}_q).
\]
But the monomial $X_1X_2^{\sqrt{q}+1}$ can't appear in
\[
g_1(\lambda X_1, X_3)X_1 + 
       g_2(\lambda X_1, X_1, X_2, X_3)X_1^2 X_2.
\]
Hence $ g_2(\lambda X_1, X_1, X_2, X_3) =0$, that is,
$X_0 - \lambda X_1$ is a factor of $g_2(X_0, X_1, X_2, X_3)$
for any $\lambda \in {\Bbb F}_q^{\ast}$.
But $\deg g_2 = \sqrt{q}-2 \,  (< q-1)$,
it is impossible.
Therefore $g_2(X_0, X_1, X_2, X_3) =0$
and $\lambda g_1(\lambda X_1, X_3)X_1 $ is Hermitian.

Let
\[
g_1(X_0, X_3) = \sum_{i=0}^{\sqrt{q}} a_i X_0^{i}X_3^{\sqrt{q}-i}.
\]
Then
\[
\lambda g_1(\lambda X_1, X_3) X_1= 
\sum_{i=0}^{\sqrt{q}} a_i \lambda^{i+1} X_1^{i+1}X_3^{\sqrt{q}-i},
\]
which is Hermitian.
Hence $a_i = 0$ for $i \neq 0, \, \sqrt{q}-1, \, \sqrt{q}$;
and also
$a_0^{\sqrt{q}} = a_{\sqrt{q}-1}$ and $a_{\sqrt{q}} \in {\Bbb F}_{\sqrt{q}}$,
that is,
\[
g_1(X_0, X_3) = a_0 X_3^{\sqrt{q}} + a_0^{\sqrt{q}}X_0^{\sqrt{q}-1}X_3 
     + a_{\sqrt{q}}X_0^{\sqrt{q}} \ 
     \mbox{\rm with $a_{\sqrt{q}} \in {\Bbb F}_{\sqrt{q}}$}.
\]
Hence
\[
F(X_0, \dots , X_3)
= a_0 X_0 X_3^{\sqrt{q}} + a_0^{\sqrt{q}}X_0^{\sqrt{q}}X_3 
     + a_{\sqrt{q}}X_0^{\sqrt{q}+1} +h(X_1, X_2, X_3), 
\]
which is Hermitian.
If $a_{\sqrt{q}} =0$,
then the surface is a cone of a nonsingular Hermitian curve
with vertex $(1,0,0,0)$, and then
$N_q(S) = (\sqrt{q}^3 +1)q +1$, which is not the given number.
So $a_{\sqrt{q}} \in {\Bbb F}_{\sqrt{q}}^{\ast}$,
and hence $S$ is a nonsingular Hermitian surface.
\qed

\vskip 22pt
\noindent
{\bf\Large Appendix}
\vskip 11pt
In this appendix, we consider a possibility of generalizing 
R\"{u}ck-Stichtenoth's characterization of the Hermitian curve
to the Hermitian surface.
\begin{atheorem}[R\"{u}ck-Stichtenoth]\label{RStheorem}
Let $C$ be a nonsingular curve of genus $g= \frac{1}{2}\sqrt{q}(\sqrt{q}-1)$
defined over ${\Bbb F}_q$.
If $N_q(C) = \sqrt{q}^3 +1$,
then the function field $\mathbb{F}_q(C)$ of $C$
coincides with that of the Hermitian curve over $\mathbb{F}_q$.
\end{atheorem}

\begin{aremark}
If one considers only nonsingular plane curves,
R\"{u}ck-Stichtenoth's characterization implies Hirschfeld et al.'s,
because the linear system $g_{\sqrt{q}+1}^2$ is unique if
$g= \frac{1}{2}\sqrt{q}(\sqrt{q}-1) > 1$.
When $q=4$, the genus in question is $1$.
As was mentioned in \cite{hir-sto-tha-vol},
the plane curve
$\{X_0 + \omega X_1 + \omega^2 X_2=0\}$
over $\mathbb{F}_4$
is not Hermitian, but has $9$ $\mathbb{F}_4$-points,
where $\mathbb{F}_4 = \{0, 1, \omega , \omega^2\}$.
However the function field of this curve over $\mathbb{F}_4$
coincides with that of the Hermitian curve over $\mathbb{F}_4$.
\end{aremark}

\ssgyokan

As has been remarked above,
the R\"{u}ck-Stichtenoth theorem gives a characterization
of the Hermitian curve in wider family than the plane curves
if one considers only nonsingular curves, namely
among the family of nonsingular curves with the same ``topological'' invariant
as the Hermitian curve.
Since the number $\sqrt{q}^3 +1$ achieves the Weil bound
for $g= \frac{1}{2}\sqrt{q}(\sqrt{q}-1) $,
the zeta function of the nonsingular curve $C$ in (A.\ref{RStheorem})
is $\frac{(1+\sqrt{q}t)^{2g}}{(1-t)(1-qt)}$.
So Theorem~A.\ref{RStheorem} can be interpreted as:
\begin{atheorem}[R\"{u}ck-Stichtenoth]
Let $C$ be a nonsingular curve of genus $g= \frac{1}{2}\sqrt{q}(\sqrt{q}-1) >1$
defined over ${\Bbb F}_q$
whose zeta function is $\frac{(1+\sqrt{q}t)^{2g}}{(1-t)(1-qt)}$.
Then $C$ is the Hermitian curve.
\end{atheorem}

On the other hand,
the number of ${\Bbb F}_q$-points of
the Hermitian surface attains the Weil-Deligne bound,
and its zeta function is:
\begin{equation}\label{zetaofhermitian}
 \frac{1}{(1-t)(1-qt)^{b_2}(1-q^2t)}  \ \ \ \ \ 
 (b_2 =q\sqrt{q} -q + \sqrt{q} + 1).
\end{equation}
It is natural to ask whether a nonsingular surface over ${\Bbb F}_q$
with the zeta function (\ref{zetaofhermitian}) is Hermitian or not.
The example below shows the answer of this question is negative.

From now on, when $X$ is a nonsingular surface over ${\Bbb F}_q$,
$Z_X(t)$ denotes the zeta function of $X$.

\begin{alemma}
Let $\Tilde{X}$ be a blowing-up of a nonsingular surface $X$ over ${\Bbb F}_q$
with center an  ${\Bbb F}_q$-point of $X$.
Then
$
Z_{\Tilde{X}}(t) = Z_X(t) \frac{1}{(1-qt)}.
$
\end{alemma}
\proof
Let $N_m(X)$ and $N_m(\Tilde{X})$ denote
the number of ${\Bbb F}_{q^m}$-points of
$X$ and $\Tilde{X}$ respectively.
Then
$N_m(\Tilde{X})= N_m(X) +q^m$.
Hence
\begin{align*}
Z_{\Tilde{X}}(t) &= \exp \left(
                     \sum_{m=1}^{\infty} (N_m(X) + q^m)\frac{t^m}{m}
                         \right)\\
                 &= \exp \left(\sum_{m=1}^{\infty}N_m(X)\frac{t^m}{m}\right) 
                      \exp\left(\sum_{m=1}^{\infty} q^m \frac{t^m}{m} \right)\\
                 &= Z_X(t) \frac{1}{(1-qt)}.
\end{align*}

\begin{aexample}
Let $X_1 = {\Bbb P}^2$ over ${\Bbb F}_q$.
It is obvious that
$
Z_{X_1}(t) = \frac{1}{(1-t)(1-qt)(1-q^2t)}.
$
Define a nonsingular surface $X_n$ over ${\Bbb F}_q$
by a blowing-up of $X_{n-1}$
with center an ${\Bbb F}_{q}$-point, successively.
Then the zeta function of $X_{b_2}$ is the same as Hermitian surface's.
However, if $\sqrt{q}+1 \geq 4$,
the Hermitian surface of degree $\sqrt{q}+1$ is not rational.
So, if $q>2^2$, then the function field of $X_{b_2}$
does not coincides with
that of the Hermitian surface over ${\Bbb F}_q$.
\end{aexample}


\end{document}